\documentclass[dvips]{amsart}
\usepackage{amssymb,amscd}
\usepackage{epsfig}
\newtheorem{theorem}{Theorem}[section]
\newtheorem{lemma}[theorem]{Lemma}
\newtheorem{corollary}[theorem]{Corollary}

\theoremstyle{definition}

\theoremstyle{remark}

\numberwithin{equation}{section}

\DeclareMathOperator{\Crit}{Crit}

\DeclareMathOperator{\End}{End}
\DeclareMathOperator{\id}{id}
\DeclareMathOperator{\const}{const}
\DeclareMathOperator{\im}{im}

\DeclareMathOperator{\e}{e}

\DeclareMathOperator{\dimcov}{\dim_{\text{\rm cov}}}
\DeclareMathOperator{\sgn}{sgn}
\DeclareMathOperator{\Hess}{Hess}
\DeclareMathOperator{\rank}{rank}

\def\bbr{\mathbb{R}}

\def\bbn{\mathbb{N}}
\def\bbz{\mathbb{Z}}

\def\Mcal{\mathcal{M}}

\def\Gcal{\mathcal{G}}
\def\Hcal{\mathcal{H}}
\def\Rcal{\mathcal{R}}

\begin{document}

\title{Equivalences for Morse homology}

\author{Matthias Schwarz}
\address{Department of Mathematics, University of Chicago,
Chicago, IL 60637}
\email{schwarz@math.uchicago.edu}
\thanks{The author was supported in part by NSF Grant \# DMS 9626430.}

\date{March 9. 1999}

\begin{abstract}
  An explicit isomorphism between Morse homology and singular homology
  is constructed via the technique of pseudo-cycles. Given a Morse
  cycle as a formal sum of critical points of a Morse function, the
  unstable manifolds for the negative gradient flow are compactified
  in a suitable way, such that gluing them appropriately leads to a
  pseudo-cycle and a well-defined integral homology class in singular
  homology.
\end{abstract}

\maketitle

\section{Introduction}
The aim of this paper is to give an explicit construction of an
isomorphism between Morse homology and singular homology.  Morse
homology is a Morse-theoretical approach to the homology of a smooth
manifold which goes back already to Thom and plays a crucial role in
Smale's proof of the $h$-cobordism theorem, cf.~also \cite{Mil-hcobordism}.
It was studied by J.~Franks \cite{Franks-MS}, rediscovered by
Witten \cite{Witten-82} in terms of a deformation of the de Rham
complex and generalized by Floer \cite{F-Witten} as an approach to
solve a conjecture by Arnold. In \cite{Sch-Morse} the author developed
a comprehensive approach to Morse homology as an axiomatic homology theory for
the category smooth manifolds (not necessarily compact) satisfying all
Eilenberg-Steenrod axioms. Moreover, this approach used the purely
relative ``Floer-theoretical'' definition of Morse homology in terms of
moduli spaces of trajectories for the gradient flow equation
connecting critical points. However, \cite{Sch-Morse} did not present an
explicit isomorphism to other axiomatic homology theories like for
instance the de Rham theorem between de Rham cohomology and singular
cohomology. That Morse homology is isomorphic to other homology
theories is proved in \cite{Sch-Morse} by extending it to a slightly
larger category of certain CW-spaces compatible with the manifold structure in which
the isomorphism is deduced inductively, based on the Eilenberg-Steenrod
axioms. That is, in such a category existence and uniqueness of the
isomorphism follows by abstract application of the axioms.

In the approach of Smale and Milnor, used similarly also in Floer's description,
a direct isomorphism between Morse homology and singular homology is
obtained by choosing a special, namely self-indexing Morse function,
such that the boundary map in the Morse chain complex can be related
to the connecting homomorphism $\partial_\ast$ in the long exact
sequence of the cell decomposition induced by the Morse function (see
also Section \ref{ssc other} below).

The approach of this paper is to show that, given {\em any} Morse function $f$
with a generic Riemannian metric $g$, one can construct singular cycles
explicitly from the given Morse cycle. The main objects which have to be
considered as intermediate tools are so-called {\em
  pseudo-cycles}. This is a geometric differential-topological way to
represent homological (integral!) cycles which plays also a role in
the definition of quantum cohomology (see e.g. \cite{McD-S}). In this
paper, we give a short proof that integral homology classes can be
represented by pseudo-cycles and that every pseudo-cycle in fact leads
to an integral homology class. 

The purpose of this paper is also to provide a detailed construction
of this equivalence between Morse homology and singular homology via
pseudo-cycles, which to the author's knowledge has not yet been
carried out elsewhere, but which already has been used several times,
in particular in the theory of quantum cohomology and Floer homology,
e.g. \cite{PSS}, \cite{Sch-qcl}, \cite{Sei}, \cite{Sch-QH}.

After a short account on the definition of Morse homology
pseudo-cycles are defined in Section \ref{sc pc}, where it is proven
that pseudo-cycles represent integral homology classes and every class
can be represented as such. Section \ref{sc expl} contains the
construction of the explicit isomorphism between Morse homology and
singular homology. In the first part we show how to obtain a
well-defined pseudo-cycle from a given Morse cycle and that the induced
singular class is uniquely associated to the Morse homology class. The
idea is to glue all unstable manifolds of critical points, which occur
in the given Morse cycle, along the $1$-codimensional strata of their
suitable compactifications. 
In the second part we construct the inverse homomorphism in terms of
intersections of pseudo-cycles representing singular classes and
stable manifolds of critical points.

\section{Morse homology}
\subsection{Definition}
Let $M$ be an oriented\footnote{If $M$ is not orientable, choose
  homology coefficients in $\bbz_2$.} smooth manifold, $f\in
C^\infty(M,\bbr)$ an ex\-haus\-ting\footnote{i.e.~proper and bounded
  below. In \cite{Sch-Morse}, this property is called coerciveness.}
Morse function and $g$ be a complete Riemannian metric.  Consider the
critical set $\Crit_\ast f$ of $f$ as graded by the Morse index
$\mu\colon \Crit f\to\bbz$ and define the stable and unstable
manifolds of the negative gradient flow in terms of spaces of curves,
\begin{equation}\label{eq unstable}\begin{split}
        W^u(x) &= \{\, \gamma\colon (-\infty,0]\to M
                  \,|\,\dot\gamma + \nabla_g f\circ\gamma=0,\,
                       \gamma(-\infty)=x  \,\},\\
        W^s(y) &= \{\, \gamma\colon [0,\infty) \to M
                  \,|\,\dot\gamma + \nabla_g f\circ\gamma=0,\,
                       \gamma(+\infty)=y  \,\}
                     \end{split}
\end{equation}
for $x,y\in\Crit f$. The curves $\gamma$ are smooth and
$\gamma(\pm\infty)$ denotes the limit for $t\to\pm\infty$.
The spaces $W^u(x)$ and $W^s(y)$ are finite-dimensional
manifolds with
$$
        \dim W^u(x)=\mu(x) \quad\text{and}\quad
        \dim W^s(y)=\dim M -\mu(y)
$$
and the evaluation mapping $\gamma\to\gamma(0)$ induces smooth
embeddings into $M$, i.e.~diffeomorphisms onto the image,
$$
        E_x\colon W^u(x)\hookrightarrow M,\quad
        E_y\colon W^s(y)\hookrightarrow M\,.
$$
However, in general, these maps are not proper.
Choosing a generic Riemannian metric we obtain Morse-Smale
trans\-ver\-sality, namely $W^u(x)$ and $W^s(y)$ intersect transversely
in $M$ with respect to $E_x$ and $E_y$. If this transversality holds
for all $x,y\in\Crit f$, $(f,g)$ is called a {\em
Morse-Smale} pair. We obtain the manifold of connecting orbits
\begin{align*}
        M_{x,y}(f,g) &=  W^u(x)\pitchfork W^s(y)\\
                     &= \{\, \gamma\colon\bbr\to M\,|\,
                        \dot\gamma+\nabla_g f\circ\gamma=0,\;
                     \gamma(-\infty)=x,\,\gamma(+\infty)=y\,
                        \},\\
        \dim M_{x,y}(f,g) &=\mu(x)-\mu(y),
\end{align*}
on which, if $x\not=y$, $\bbr$ acts freely and properly by shifting
$$
        (\tau\ast\gamma)(t)=\gamma(t+\tau)\,.
$$

Let us fix orientations for all unstable manifolds $W^u(x)$, then the
orientation of $M$ induces orientations for $W^s(y)$ and $M_{x,y}$. We
call an unparameterized trajectory $\hat\gamma\in M_{x,y}/\bbr$ for
relative index $1$ positively oriented if the orbit
$\bbr\cdot\hat\gamma\subset M_{x,y}$ is positively oriented by the
action of $\bbr$ which corresponds to the action by the negative
gradient flow. 
Thus, for relative index $1$, the moduli spaces of
unparameterized trajectories
$$
        \widehat M_{x,y}=M_{x,y}/\bbr,\quad \mu(x)-\mu(y)=1,
$$
are compact, that is finite, and every element $\hat\gamma$ carries a
sign $\tau(\hat\gamma)\in\{\pm1\}$. We define the intersection numbers
$$
        n(x,y) = \sum_{\hat\gamma\in\widehat M_{x,y}}
                 \tau(\hat\gamma)
$$
and an operator on the module over $\bbz$ generated by the critical points of
index $k$, 
\begin{gather*}
        C_k(f) = \bbz\otimes \Crit_k f,\quad
        \partial=\partial(f,g),\\
        \partial\colon C_k(f) \to C_{k-1}(f),\quad
        \partial x = \sum_{y} n(x,y)\,y\,.
\end{gather*}
The fundamental theorem of Morse homology is
\begin{theorem}\label{thm boundary}
        $\partial$ is a chain boundary operator, i.e.~$\partial\circ\partial=0$.
\end{theorem}
Hence, the homology $H_k(f,g;\bbz)=H_k(C_\ast(f),\partial(f,g))$ is
well-defined as the quotient of the module of {\em Morse-cycles}
$$
        Z_k(f,g) = \{\, a=\hspace{-1em}\sum_{x\in\Crit_k f}\hspace{-1em} a_x x \,|\,
        \partial a=0\,\}\,.
$$
modulo the boundaries $B_k(f,g)=\im\partial$.

Let us now recall the homotopy invariance result in Morse homology. It
is based on Conley's continuation principle (see \cite{Con-isolated}).
\begin{theorem}\label{thm continuation}
  Given two Morse-Smale pairs $(f^0,g^0)$ and $(f^1,g^1)$ there exists a
  canonical homomorphism
$$
        \Phi_{10}\colon H_\ast(f^0,g^0) \to H_\ast(f^1,g^1)
$$
  such that
$$
        \Phi_{21}\circ\Phi_{10} = \Phi_{20}\quad
        \text{and}\quad  \Phi_{00}=\id\,.
$$
  In particular, every $\Phi_{ji}$ is an isomorphism.
\end{theorem}
This continuation theorem implies that we have well-defined Morse homology groups
\begin{equation}\label{eq Morse homology}\begin{split}
        &H^{\text{Morse}}_\ast(M;\bbz) \stackrel{\mathsf{def}}{=} \big\{\,
        (a_i)\in\prod H_\ast(f^i,g^i)\,|\,
        a_j=\Phi_{ji}a_i\,\big\}\\
        &\rho_i\colon
        H_\ast(f^i,g^i)\stackrel{\cong}{\longrightarrow}H^{\text{Morse}}_\ast(M;\bbz),\quad
        \rho_i\circ\Phi_{ij}=\rho_j\,.
      \end{split}
\end{equation}
Let us recall the construction of $\Phi_{10}$ from \cite{Sch-Morse}.
Given the Morse-Smale pairs $(f^i,g^i)$, $i=0,1$, we choose an
asymptotically constant homotopy over $\bbr$, $(f_s,g_s)$, $s\in\bbr$
with
$$
        (f_s,g_s) = \begin{cases}
        (f^0,g^0) ,& s\leqslant -R,\\
        (f^1,g^1) ,& s\geqslant R,
                    \end{cases}
$$
for $R$ large enough. This gives rise to the trajectory spaces
\begin{align*}
        M_{x_o,x_1}(f_s,g_s) = \{\,\gamma\,|\, &\dot\gamma(s)
        + \nabla_{g_s}f_s(\gamma(s))=0,\\
        &\gamma(-\infty)=x_o,\,\gamma(\infty)=x_1\,\}\,.
\end{align*}
For a generic choice of the homotopy $(f_s,g_s)$, these spaces
are finite dimensional manifolds with
$$
        \dim M_{x_o,x_1}=\mu(x_o)-\mu(x_1)
$$
and compact in dimension $0$. As in the definition of the boundary
operator $\partial$ we define
\begin{gather*}
        \Phi_{10}\colon C_\ast(f^0,g^0) \to C_\ast(f^1,g^1),\\
        \Phi_{10} x_o = \sum_{x_1} n(x_o,x_1) x_1,\\
        n(x_o,x_1)=\#_{\text{alg}}M_{x_o,x_1}(f_s,g_s)\,,
\end{gather*}
where $\#_{\text{alg}}$ means counting with signs $\tau(u)=\pm1$,
analogously to above. That $\Phi_{10}$ is well-defined on the level of
homology follows from a theorem stating that
$$
        \Phi_{10}\circ\partial_0 = \partial_1\circ\Phi_{10}\,.
$$
Moreover, it is shown in \cite{Sch-Morse} that the homomorphism
$\Phi_{10}$ on homology level does not depend on the choice of the
homotopy $(f_s,g_s)$.

\subsection{Isomorphism via axiomatic approach}\label{ssc axiom}
In \cite{Sch-Morse}, Morse homology is extended towards an
axiomatic homology theory for the category of smooth manifolds. It is
functorial with respect to smooth maps, there exists a relative
version so that we have an associated long exact sequence, and all
axioms of Eilenberg and Steenrod are satisfied. However, in order to
derive a natural isomorphism with any other axiomatic homology theory,
an extension to a larger category of spaces is required, e.g.~towards
the subcategory of CW-pairs which are embedded smoothly into
finite-dimensional manifolds as strong deformation retracts of open
subsets. This approach is adopted in \cite{Sch-Morse} in order to
prove the equivalence with other homology theories.

\section{Pseudo-cycle homology}\label{sc pc}
In \cite{McD-S}, pseudo-cycles were defined in order to find a suitable
dif\-fe\-ren\-tial-to\-po\-lo\-gi\-cal representation of homology
cycles. However, this was only used with rational coefficients so that
every cycle can be represented as a closed submanifold. Here, we
consider integral homology classes.

Let $M$ be a compact\footnote{This poses no restriction for our
  application to Morse homology because we consider only cycles lying
  in the compact sublevel sets $M^a=\{\,p\in M\,|\,f(p)\leqslant a\,\}$ of 
  an exhausting function.} $n$-dimensional manifold. We consider an oriented
smooth $k$-dimensional manifold without boundary $V$ together with a
smooth map $f\colon V\to M$. Let the set $f(V^\infty)$ be defined as
in \cite{McD-S},
\begin{equation}\label{eq pseudocycle}
        f(V^\infty)\stackrel{\mathsf{def}}{=} \bigcap_{K\subset V\text{ cpt.}}
        \overline{f(V\setminus K)}\,.
\end{equation}
According to \cite{McD-S}, $f\colon V \to M$ is a {\bf pseudo-cycle} if
$f(V^\infty)$ can be covered by the image of a smooth map $g\colon P \to M$
which is defined on a manifold $P$ of dimension not larger than $\dim
V-2$. 

Moreover, let $W$ be an oriented smooth $(k+1)$-dimensional manifold
with boundary $\partial W$, such that the inclusion $i\colon \partial
W\hookrightarrow W$ is proper, and let $F\colon W\to M$ be a smooth map.
\begin{theorem}\label{prop cycle1}
Let $(f,V)$ be a pseudo-cycle and $(F,W)$ as above.
\begin{itemize}
        \item[(a)] 
        If $H_k(f(V^\infty);\bbz)=H_{k-1}(f(V^\infty);\bbz)=0$ and
        $f(V)\not\subseteq f(V^\infty)$, then $(f,V)$ induces a unique
        integral homology class $\alpha_f\in H_k(M;\bbz)$.
        \item[(b)]
        Let $\partial W=U$ be an open subset of $V$ such that
        $f(V^\infty)\subseteq f(U^\infty)$ and $F(W^\infty)\cap
        f(U)=\emptyset$. If $H_k(F(W^\infty);\bbz)=0$ the
        homology class $\alpha_f$ vanishes.
\end{itemize}
\end{theorem}
In view of part (b) let us consider two pseudo-cycles $f_1\colon V_1\to
M$ and $f_2\colon V_2\to M$ to be {\bf cobordant} if their disjoint
union $V=V_1\amalg V_2^\ast$ with orientation on $V_2$ reversed forms
a pseudo-cycle $f\colon V\to M$ such that there exists $F\colon W\to M$ 
satisfying the condition in (b).

In this section we are using Alexander-Spanier homology theory for
locally compact Hausdorff spaces, cf.~\cite{Massey}. The homology
theory with arbitrary supports\footnote{Alternatively, we could also
  use Borel-Moore homology with specified type of supports.}, i.e.~not
necessarily compact supports, is denoted by $H^\infty_\ast(X)$.  Note
that, however, homology theory with arbitrary supports, which is
functorial with respect to proper maps, agrees with any homology
theory with compact supports, as for instance singular homology, when
restricted to compact sets as $M$, $f(V^\infty)$ and $F(W^\infty)$.

\begin{proof}[Proof of Theorem \ref{prop cycle1}]
Every oriented $k$-dimensional manifold $X$ without boundary carries a
uniquely defined fundamental class $[X]\in H^\infty_k(X;\bbz)$. If $X$
is a manifold with boundary $\partial X$ then $[X]$ is well-defined in
$H^\infty_k(X\setminus\partial X)=H_k^\infty(X,\partial X)$.  Every open
subset $U\subset X$ inherits an orientation from $X$ so that the natural
restriction map $\rho\colon H^\infty_k(M;\bbz) \to H^\infty_k(U;\bbz)$
gives $\rho([U])=[X]$. Without loss of generality we may assume that
$$
        f(V)\cap f(V^\infty) =\emptyset\,.
$$
Otherwise, we replace $V$ by the open, nonempty subset $V\setminus
f^{-1}(f(V^\infty))$. Since Alexander-Spanier homology with arbitrary
supports is functorial with respect to proper maps of locally compact
Hausdorff spaces we redefine the map
$$
        f\colon V \to M\setminus f(V^\infty)\,.
$$ 
By definition of $f(V^\infty)$, $f$ is proper. The integral class
$$
        (f)_\ast([V])\in H^\infty_k(M\setminus f(V^\infty);\bbz)
$$
is well-defined. From the exact homology sequence for
$H^\infty_\ast$ and the pair $(M,f(V^\infty))$,
$$
        H_k(f(V^\infty)) \to
        H_k(M) \stackrel{j_\ast}{\longrightarrow}
        H^\infty_k(M\setminus f(V^\infty)) \to
        H_{k-1}(f(V^\infty))\,,
$$
we obtain by assumption the isomorphism $j_\ast$. Hence,
$$
        \alpha_f\equiv j_\ast^{-1}(f)_\ast([V])
        \in H_k(M;\bbz)
$$
is well-defined. 

We consider now an open subset $U\subset V$ such that
$f(V^\infty)\subset f(U^\infty)$ and $f(U)\cap
f(U^\infty)=\emptyset$. We carry out the same procedure as before for
the proper map
$$
        f_U\colon U \to M\setminus f(U^\infty)
$$
and relate it to $\alpha_f$ by the following commutative diagram with
respect to the natural restriction homomorphism $\rho$,
$$
\begin{CD}
        H^\infty_k(V) @>{f_\ast}>> H^\infty_k(M\setminus f(V^\infty))\\
        @VV{\rho}V                 @VV{\rho}V\\
        H^\infty_k(U) @>{(f_U)_\ast}>>
        H^\infty_k(M\setminus f(U^\infty))\,.
\end{CD}
$$
Since $\rho\circ j_\ast=j^U_\ast$ it follows that
\begin{equation}\label{eq commut}
        j^U_\ast(\alpha_f) = \rho\circ f_\ast([V])
        = (f_U)_\ast([U])\,.
\end{equation}
                
Let us consider now the bordism $\partial W=U$. Without loss of
generality we can assume that $F(W^\infty)\cap F(W)=\emptyset$ so that
we have the proper map 
$$
        F\colon W\to M\setminus F(W^\infty)\,.
$$
In Alexander-Spanier homology theory, we know that the fundamental class
$[U]$ is the image of the fundamental class of the manifold with
boundary $W$ under the boundary homomorphism $\partial_\ast$ in the
exact homology sequence of the pair $(W,U)$. We obtain the
commutative diagram
$$
\begin{CD}
        H_{k+1}^\infty(W,U) @>{\partial_\ast}>>
        H_k^\infty(U) @>{i_\ast}>> H_k^\infty(W)\\
        && @VV{(f_U)_\ast}V     @VV{F_\ast}V\\
        && H^\infty_k(M\setminus f(U^\infty))  @>{\rho}>>
        H^\infty_k(M\setminus F(W^\infty))\\
        && @AA{j^U_\ast}A               @AA{j^W_\ast}A\\
        && H_k(M)       @>{\id}>>       H_k(M)\,,
\end{CD}
$$
therefore by (\ref{eq commut})
$$
        j^W_\ast(\alpha_f)=\rho\circ j^U_\ast(\alpha_f)
        =F_\ast\circ i_\ast([U])=0,
$$
because $[U]\in\im\partial_\ast$. Since $j^W_\ast$ is injective due to
$H^\infty_k(F(W^\infty))=0$ it follows that $\alpha_f$ vanishes.
\end{proof}

A topological space $S$ is said to have {\bf covering dimension} at most
$n$ if every open cover $\mathfrak{U}=\{U_\alpha\}$ has a refinement
$\mathfrak{U}'=\{U_\alpha'\}$ for which all the $(n+2)$-fold
intersections are empty\footnote{compare \cite{Don-Kron}, Section 9.2.3.},
$$
        U'_B = \bigcap _{\beta\in B} U'_\beta = \emptyset
        \quad\text{if } |B| \geqslant n+2\,.
$$
We say then ${\dim_{\text{\rm cov}}} S \leqslant n$. Clearly, $S$ having covering
dimension at most $n$ implies that $\check H^m(S)=0$ for $m>n$. For a
compact space $S$ this implies $H_m(S)=0$ for $m>n$.

We have the following simple
\begin{lemma}\label{lm covering}
  Let $f\colon P\to M$ be a smooth map between manifolds and $S$ be a
  compact subset of $M$ such that $S\subset f(P)$. Then
$$
        \dimcov S \leqslant \dim P\,.
$$
\end{lemma}
The final result is
\begin{theorem}\label{cor pseudo-cycle}
  Every pseudo-cycle $f\colon V\to M$ of dimension $k$ induces a
  well-defined integral homology class $\alpha_f\in
  H_k(M;\bbz)$. Moreover, any singular cycle $\alpha\in
  Z^{\mathsf{sing}}_k(M;\bbz)$ gives rise to a $k$-pseudo-cycle
  $f\colon V\to M$ such that $\alpha_f=\alpha$.
\end{theorem}
Note that if $f(V)\subset f(V^\infty)$ then trivially $\alpha_f=0$.
\begin{proof}
The well-defined homology class $\alpha_f$ follows from combining
Theorem \ref{prop cycle1} with Lemma \ref{lm covering}. 

Suppose now that $\alpha\in H_k^{\mathsf{sing}}(M;\bbz)$ is a
$k$-cycle given by a smooth singular chain. By pairwise identifying
and sufficiently smoothing the
$k-1$-dimensional faces of the $k$-simplexes involved in $\alpha$ the
cycle-property of $\alpha$ implies that we obtain a $k$-dimensional manifold $V$, not
necessarily compact, with a smooth structure, such that the singular
chain gives a map $f\colon V\to M$ meeting the pseudo-cycle
condition, since $f(V^\infty)$ is covered by the images of the faces
of codimension $2$ and higher. Two cohomologous singular chains lead
to cobordant pseudo-cycles in the sense of Theorem \ref{prop cycle1} (b).
\end{proof}

Hence, from now on we can represent integral cycles in singular
homology by pseudo-cycles.

\section{The Explicit Isomorphism}\label{sc expl}
The first part consists of showing that each Morse cycle leads to a
well-defined pseudo-cycle, and that the associated singular homology
class does not depend on any uncanonical choices involved.
\subsection{Pseudo-cycles from Morse cycles}
Let $(f,g)$ be a Morse-Smale pair and consider the associated homology
$H_\ast(f,g)$. The idea of defining the homomorphism into singular
homology is to construct a $k$-dimensional pseudo-cycle $E\colon
Z(a)\to M$ for a given Morse cycle $\{a\}=\{\sum_{x\in\Crit_k f}
a_x\,x\in\}\in H_k(f,g)$. This is essentially based on considering the 
unstable manifolds $W^u(x)$ from (\ref{eq unstable}) with multiplicity $a_x\in\bbz$ and their
evaluation maps $E_x\colon W^u(x)\to M$. In order to obtain a
well-defined pseudo-cycle we have to carry out a suitable
identification on the $1$-codimensional strata of a suitable
compactification of $W^u(x)$.

Let $x\in\Crit_k f$ and $y\in\Crit_{k-1}f$ such that $\widehat M_{x,y}$
is a nonempty finite set. We say that a sequence $(w_n)\subset W^u(x)$
is {\em weakly convergent} towards a simply broken trajectory,
$$
        w_n\rightharpoonup (\hat u,v)\in \widehat M_{x,y}\times W^u(y),
$$
if $w_n\to v$ in $C^\infty_{\text{loc}}((-\infty,0],M)$ and there
exists a reparametrization sequence $\tau_n\to-\infty$ such that
$\tau_n\ast w_n\to u$ in $C^\infty$ on compact subsets of $\bbr$ for a
representative $u$ of the unparameterized trajectory $\hat u$. Note that,
in particular, $w_n(0)\to v(0)$.
\begin{figure}[htb] 
 \begin{center}
\setlength{\unitlength}{0.1mm}%
\begin{picture}(600,600)(0,0)
 \put(0,0){\epsfig{file=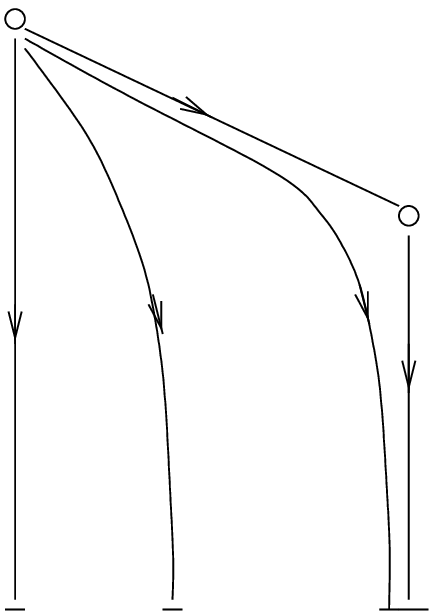,height=6cm}}
 \put(50,590){\makebox(0,0)[lb]{$x$}}
 \put(420,375){\makebox(0,0)[lb]{$y$}}
 \put(250,500){\makebox(0,0)[lb]{$\hat u$}}
 \put(415,180){\makebox(0,0)[lb]{$v$}}
 \put(90,300){\makebox(0,0)[lb]{$w_n$}}
 \put(70,0){\makebox(0,0)[lb]{$w_n(0)$}}
 \put(420,0){\makebox(0,0)[lb]{$v(0)$}}
\end{picture}
 \end{center}
 \caption{Weak convergence towards simply broken trajectory
\label{fig simply-broken}}
\end{figure}

The following result is completely analogous to the gluing results
developed in \cite{Sch-Morse}, Section 2.5. There, the gluing operation has been
constructed for trajectory spaces $M_{y,z}$ instead of $W^u(y)$, but
the case of unstable manifolds is handled exactly the same. It
provides us with the suitable description of strata of the weak
compactification of $W^u(x)$.

\begin{lemma}[\cite{Sch-Morse}]\label{prop gluing1}
  Given an open subset $V\subset W^u(y)$ with compact closure
  there exists a constant $\rho_V>0$ and a smooth map
\begin{align*}
        \#^V \colon &\widehat M_{x,y}\times V\times
        [\rho_V,\infty) \to W^u(x),\\
        &(\hat u,v,\rho) \mapsto \hat u\#_\rho v,
\end{align*}
  such that
\begin{itemize}
\item[(a)] $\#^V$ is an embedding,
\item[(b)] $\#^V(\hat u,\cdot,\cdot)$ is orientation preserving
        exactly if $\tau(\hat u)=+1$,
\item[(c)] $\hat u\#_\rho v\rightharpoonup (\hat u,v)$ for
        $\rho\to\infty$, and for any $w_n\rightharpoonup(\hat u,v)$ there
        exists an $n_o$ such that for all $n\geqslant n_o$ $w_n=\hat
        u_n\#_{\rho_n} v_n$ for unique $(\hat u_n,v_n,\rho_n)$, and
\item[(d)] the evaluation maps $E_x\colon W^u(x)\to M$ and
        $E_y\colon W^u(y)\to M$ extend to
$$
        \bar E_x\colon W^u(x) \cup_{\#^V}
        \widehat M_{x,y}\times V\times [\rho_V,\infty)\to M\,.
$$
               such that $\bar E_x(\hat u,v,\rho)=E_x(\hat u\#_\rho
               v)$ for $\rho\in[\rho_v,\infty)$ and $\bar E_x(\hat
               u,v,\rho)=E_y(v)$ for $\rho=\infty$.
\end{itemize}
\end{lemma}
Let us define $\overline{W}{}^u(x)$ to be the disjoint union
\begin{equation}\label{eq boundary}
        \overline{W}{}^u(x) = W^u(x) \cup \bigcup_{\mu(y)=\mu(x)-1}
        \widehat M_{x,y}\times W^u(y)
\end{equation}
equipped with the topology generated by
\begin{itemize}
\item[(a)] the open subsets of $W^u(x)$,
\item[(b)] the neighborhoods of $(\hat u,v)\in\widehat M_{x,y}\times
        W^u(y)$ of the form 
$$
        \#^V(\{\hat u\}\times V\times (\rho,\infty)) \cup \{\hat
        u\}\times V,\quad \rho\geqslant\rho_V,
$$ 
        for $V\subset W^u(y)$ open with compact closure.
\end{itemize}
This provides a Hausdorff topology and we obtain
\begin{lemma}\label{prop gluing2}
        $\overline{W}{}^u(x)$ is an oriented manifold with boundary oriented
        by $\widehat M_{x,y}\times W^u(y)$ and $\bar
        E_x\colon\overline{W}{}^u(x)\to M$ is a smooth embedding.
\end{lemma}\noindent
The proof is given below together with the proof of Lemma \ref{prop cycle}.

Consider now a Morse-cycle 
$$
        a\in Z_k(f,g),\quad a=\sum_x a_xx\,.
$$
Given $l\in\bbn$ let us denote by $l\cdot\overline{W}{}^u(x)$ the disjoint
union of $l$ copies of $\overline{W}{}^u(x)$, that is, the topological
sum. If $l\in\bbz$, $l<0$, we replace $\overline{W}{}^u(x)$ by
$\overline{W}{}^u(x)^\ast$, that is, with the orientation reversed. Thus,
we associate to $a$ the topological sum
\begin{equation}\label{eq topsum}
        a\mapsto \amalg_x a_x\cdot\overline{W}{}^u(x)
\end{equation}
which is a $k$-dimensional oriented manifold with oriented boundary
and it consists of $\sum_x |a_x|$ connected components. Observe that
this manifold with boundary is not compact in general.

We denote by $\Delta a$ the following finite set of connecting
un\-para\-me\-trized trajectories of relative index $1$,
$$
        \Delta a = \bigcup \{\, a_x\hat u\,|\,
        \hat u\in\widehat M_{x,y},\,x\in\Crit_k f,\,
        y\in\Crit_{k-1}f\,\}
$$
where $a_x\hat u$ is the disjoint union of $|a_x|$ copies of $\{\hat
u\}$. Each $\hat u$ carries the sign $\tau(\hat u)\in\{\pm1\}$ and we
assign to every $\gamma\in a_x\hat u$ the new sign $\sigma(\gamma)=\sgn(a_x)\cdot
\tau(\hat u)$. Computing
$$
        \partial a = \sum_x \sum_{\mu(y)=\mu(x)-1}
                \sum_{\hat u\in\widehat M_{x,y}}
                a_x\tau(\hat u)\,y
$$
we immediately obtain
\begin{lemma}\label{lm 1-1}
  If $a=\sum_x a_x x$ is a Morse-cycle there exists an equivalence
  relation $\sim_{\Delta a}$ on $\Delta a$ such that for each
  $\gamma\in\Delta a$ there exists a unique $\gamma'\not=\gamma$ with
  $\sigma(\gamma')=-\sigma(\gamma)$, so that $\gamma\sim_{\Delta a}\gamma'$ and
  $\gamma(+\infty)=\gamma'(+\infty)\in\Crit_{k-1}f$ for $\gamma,\gamma'$ viewed as flow trajectories.
\end{lemma}
Since $\Delta a$ is an index set for the components of the
$k-1$-dimensional manifold from (\ref{eq topsum}),
$\partial(\amalg_x a_x\overline{W}{}^u(x))$ such that $\sigma(\gamma)$
corresponds to the boundary orientation, we obtain the equivalence
relation for points $\{\gamma\}\times\{v\}\in a_x \widehat M_{x,y}\times
W^u(y)$,
$$
        \{\gamma\}\times \{v\} \sim_a \{\gamma'\}\times\{v'\}
        \quad\stackrel{\mathsf{def}}{\Longleftrightarrow}\quad
        \gamma \sim_{\Delta a} \gamma',\; v=v'\,.
$$
We define
\begin{equation}\label{eq cycle}
        Z(a) = \amalg_x a_x\overline{W}{}^u(x) \big/ \sim_a\,.
\end{equation}
\begin{figure}[htb] 
 \begin{center}
\setlength{\unitlength}{0.1mm}%
\begin{picture}(900,500)(0,0)
 \put(0,0){\epsfig{file=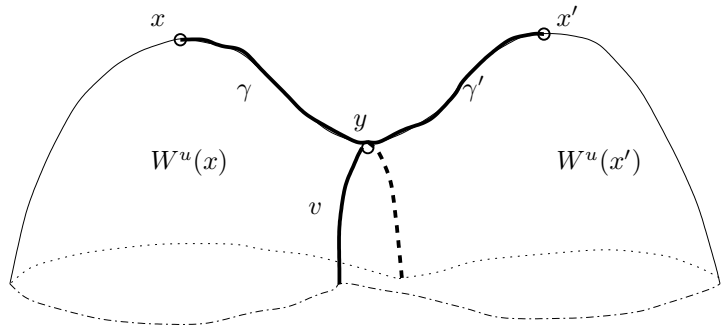,height=4cm}}
 \put(190,400){\makebox(0,0)[lb]{$x$}}
 \put(190,200){\makebox(0,0)[lb]{$W^u(x)$}}
 \put(730,400){\makebox(0,0)[lb]{$x'$}}
 \put(730,200){\makebox(0,0)[lb]{$W^u(x')$}}
 \put(460,260){\makebox(0,0)[lb]{$y$}}
 \put(305,300){\makebox(0,0)[lb]{$\gamma$}}
 \put(605,300){\makebox(0,0)[lb]{$\gamma'$}}
 \put(400,150){\makebox(0,0)[lb]{$v$}}
\end{picture}
 \end{center}
 \caption{Gluing unstable manifold along simply broken trajectories
\label{fig cycle}}
\end{figure}
One easily sees that $Z(a)$ is a topological Hausdorff space and
clearly the evaluation maps $\bar E_x$ yields
$$
        E\colon Z(a)\to M,\quad
        [\gamma,v] \mapsto v(0)\,.
$$
In fact, we obtain
\begin{lemma}\label{prop cycle}
  The space $Z(a)$ carries the structure of a $k$-dimen\-sio\-nal
  manifold without boundary and $E\colon Z(a)\to M$ is a smooth map.
\end{lemma}

\begin{proof}[Proof of Lemmata \ref{prop gluing2} and \ref{prop cycle}]
Let us consider $x,x'\in\Crit_k f$ and $y\in\Crit_{k-1}f$ with $\hat
u\in\widehat M_{x,y}$ and $\hat u'\in\widehat M_{x',y}$ such that $\hat
u\sim_{\Delta a}\hat u'$, in particular, $\tau(\hat u)=-\tau(\hat u')$. Let
$v_o\in W^u(y)$ and $V,V'\subset W^u(y)$ be two relatively compact
neighborhoods of $v_o$.

In view of Lemma \ref{prop gluing2} we consider the following
local coordinates at $(\hat u,v_o)\in\partial\overline{W}{}^u(x)$,
respectively for Lemma \ref{prop cycle} $[(\hat
u,v_o)]_{\sim_a}\in Z(a)=\amalg_x a_x\overline{W}{}^u(x)/\sim_a$,
\begin{gather*}
        \#^V_{\hat u,\hat u'}\colon V\times [0,\epsilon) \to \overline
        {W}^u(x),\\
        (v,t) \mapsto \begin{cases}
        \hat u \#_{-\frac{1}{t}} v, & t<0,\\
        (\hat u,v_o) \sim_a (\hat u',v_o), &t=0,\\
        \hat u' \#_{\frac{1}{t}} v, & t>0,
                      \end{cases}
\end{gather*}
where $\#$ is the gluing map from Proposition \ref{prop gluing1} and
$\epsilon>0$ is small enough depending on the compact set
$\operatorname{cl}(V)$. Thus, we have to show that
\begin{itemize}
\item[(A)]
        $(\#^{V'}_{\hat u,\hat u'}){}^{-1}\circ \#^V_{\hat u,\hat u'}
        \colon U \times (-\epsilon_o,\epsilon_o) \to (V\cap V')
        \times \bbr$
        is smooth for $U\subset V\cap V'$ and
        $\epsilon_o<\min(\epsilon,\epsilon')$ sufficiently small, and
        that  
\item[(B)]
        $E\circ\#^V_{\hat u,\hat u'}\colon V\times
        (-\epsilon,\epsilon)\to M$ is smooth at $(v,0)$.
\end{itemize}

Let us recall the definition of $\hat u\#_\rho v$ from
\cite{Sch-Morse}. Let $\beta^-\colon\bbr\to[0,1]$ be a cut-off
function with
$$
        \beta^-(s)=\begin{cases}
        1, & s\leqslant -1,\\
        0, & s\geqslant 0,
                   \end{cases}
$$
and $\beta^+(s)=\beta^-(-s)$. We write
$$
        \beta^\pm_\rho(s)=\beta^\pm(s+\rho), \quad
        \hat u_\rho(s)=\hat u(s+\rho)\,.
$$
For every $v\in V$ and $\rho\geqslant \rho_o$ large enough we define
$w=w(\hat u,v,\rho)$ by
$$
        w=\begin{cases}
        \hat u(s+2\rho),        &s\leqslant -\rho-1,\\
        \exp_y\big(  \beta^-_\rho \exp_y^{-1}\circ\hat u_{2\rho}
                    +\beta^+_\rho \exp_y^{-1}\circ v \big)(s),
                                & |s+\rho|<1,\\
        v(s),                   &s\geqslant -\rho+1\,.
          \end{cases}
$$
In particular, $w(\rho)=y$. One can find a $\rho_V>0$ and a bundle
$\pi\colon L^\bot\to V\times [\rho_V,\infty)$ with $L_{(v,\rho)}^\bot
\subset C^\infty(w^\ast TM)$ such that there exists a unique
section $\gamma\colon V\times [\rho_V,\infty) \to L^\bot$ providing
$$
        (\hat u\#_\rho v)(s) =\exp_{w(s)} (\gamma(v,\rho)(s)),
        \quad \hat u\#_\rho v\in W^u(x)\,.
$$
The bundle $L^\bot$ can be completed fiberwise in terms of a Sobolev
space yielding a smooth bundle such that $\gamma$ is a smooth
section. (Details can be found in \cite{Sch-Morse}.)
Moreover, there is an exponential estimate for the correction term
$\gamma(v,\rho)$ between $w(v,\rho)$ and $\hat u\#_\rho v$. Namely, there
exists a $\sigma>0$ such that
\begin{equation}\label{eq exponential1}
        \sup_{s\in\bbr} |\gamma(v,\rho)(s)| \leqslant
        c \e^{-\sigma\rho}
\end{equation}
for some $c>0$ uniformly for $v\in V$. Moreover, also the covariant derivatives
of $\gamma$ with respect to $v$ and $\rho$ satisfy such an exponential
estimate as $\rho\to\infty$. This is due to the fact that
$w(v,\rho)\rightharpoonup (\hat u,v)$ as $\rho\to\infty$ and that the
gradient flow trajectories $\hat u$ and $v$ converge exponentially fast
towards $y$,
$$
        d(\hat u(s),y),\, d(v(-s),y) \leqslant
        c\e^{-\sigma s} \quad\text{as }s\to\infty\,.
$$
The construction of $L^\bot$ and $\gamma$ in \cite{Sch-Morse} is such
that $L^\bot\to V\times [\rho_V,\infty)$ and $L^\bot\to V'\times
[\rho_{V'},\infty)$ coincide over $V\cap V'\times
[\max(\rho_V,\rho_{V'}),\infty)$. One obtains a unique smooth gluing
map
$$
        \#^{V\cup V'}_{\hat u}\colon (V\cup V')\times
        [\max(\rho_V,\rho_{V'}),\infty)\to W^u(x)
$$
extending $\#^V$ and $\#^{V'}$ and assertion (A) follows.

Let us consider now the coordinate chart $\phi_V(v,t)=\#^V_{\hat
u,\hat u'}(v,t)$ with the exponential estimate for the correction term
$\gamma(v,\pm\frac{1}{t})$,
\begin{equation}\label{eq exponential2}\textstyle
        \| \nabla^\alpha\gamma(v,\pm\frac{1}{t})\|_\infty \leqslant
        c_\alpha \e^{-|\frac{1}{t}|}\,
\end{equation}
where $\nabla^\alpha$ are the covariant derivatives of the section
$\gamma$ with respect to the variables $v$ and $\rho$.
We obtain for the evaluation map $E\colon Z(a)\to M$,
\begin{align*}
        E \circ \phi_V(v,t)     &=\begin{cases}
                E_x(\hat u'\#_{(-\frac{1}{t})} v), & t<0,\\
                E_y(v),                            & t=0,\\
                E_x(\hat u \#_{\frac{1}{t}} v),    & t>0,
                                  \end{cases}\\
        &= \begin{cases}
        \exp_{v(0)} \big( \gamma(\hat u',v,-\frac{1}{t})(0)\big), 
              & t<0,\\
        v(0), & t=0,\\
        \exp_{v(0)} \big( \gamma(\hat u,v,\frac{1}{t})(0) \big),
              & t>0\,.
          \end{cases}
\end{align*}
Thus, the smoothness of $E\circ\phi_V$ follows from 
(\ref{eq exponential2}) 
and the standard identities for the covariant
derivatives of $\exp\colon TM\to M$ at $0_p\in T_p M$.
\end{proof}

The next step is to analyze the map $E\colon Z(a) \to M$ with respect to
the end of $Z(a)$, because in general $Z(a)$ is not
compact. If it is compact, we immediately obtain the well-defined
integral homology class $E_\ast([Z(a)])\in H_k(M;\bbz)$ associated to
the Morse cycle $a$ of degree $k$. 
\begin{lemma}\label{prop pseudo-cycle}
  The evaluation map $E\colon Z(a) \to M$ associated to a Morse cycle
  $a\in Z_k(f,g)$ is a $k$-dimensional pseudo-cycle.
\end{lemma}

\begin{proof}
Consider a point $p\in M$ such that
$$
        p\in E(Z(a)^\infty)=\bigcap_{K\underset{\text{cpt.}}{\subset}Z(a)}
        \overline{ E(Z(a)\setminus K)}\,.
$$
That is, there exists a sequence $(\gamma_n)\subset Z(a)$ such that
$\gamma_n(0)\to p$ in $M$ but $(\gamma_n)$ contains no convergent
subsequence in $Z(a)$. We can assume that every $\gamma_n$ corresponds to an
element in $W^u(x)$ for some $x\in \Crit_kf$ such that $a_x\not=0$ for
$a=\sum_x a_x x$. The compactness result for the space of negative
gradient flow trajectories provides a convergent subsequence
$$
        \gamma_{n_k} \stackrel{C^\infty_{\text{loc}}}
        {\longrightarrow} \gamma \in W^u(z),
$$
with $\mu(z)\leqslant \mu(x)$. Since $(\gamma_{n_k})$ does not converge
in 
$$
        Z(a)=\amalg_x a_x \overline{W}{}^u(x) \big/\sim\,,
$$
we obtain $\mu(z) \leqslant \mu(x) -2$. This shows that
$E(Z(a)^\infty)$ is covered by the images of the evaluation maps
\begin{equation}\label{eq cover}
        E(Z(a)^\infty) \subset
        \bigcup\limits_{\mu(z)\leqslant k-2}
        \im E_z \subset M\,.
\end{equation}
Thus, $E\colon Z(a)\to M$ is a pseudo-cycle.
\end{proof}

Given a Morse-cycle $a\in Z_k(f,g)$, we denote the by Lemma \ref{prop
  pseudo-cycle} and Theorem \ref{cor 
  pseudo-cycle} uniquely determined homology class by
$$
        [a]\in H_k(M;\bbz)\,.
$$
Moreover, the map $a\mapsto [a]$ is linear by construction.
\begin{theorem}\label{prop cycle2}
  If $a\in Z_k(f,g)$ is a boundary, i.e.~$a=\partial b$ for some $b\in
  C_{k+1}(f,g)$, then $[a]=0$. That is, the homomorphism
\begin{equation}\label{eq homom}
        \Phi_{f,g}\colon H_k(C_\ast(f,g),\partial) \to H_k(M;\bbz),\quad
        \{a\} \mapsto [a]
\end{equation}
  is well-defined.
\end{theorem}
\begin{proof}
Consider $a=\sum_x a_x x\in Z_k(f,g)$ and $b=\sum_z b_z z\in
C_{k+1}(f,g)$ such that $\partial b=a$. Similar to the construction of
$Z(a)$ we now set
$$
        W = \amalg_{z\in \Crit_{k+1}f}
        b_z\cdot \overline{W}{}^u(z), \quad
        \overline{W}{}^u(z)= W^u(z) \cup \bigcup\limits_{x\in\Crit_k f}
        \widehat M_{z,x}\times W^u(x)\,,
$$
as in (\ref{eq boundary}). We can obtain the boundary of the manifold
$\overline{W}{}^u(z)$ as
$$
        \partial\overline{W}{}^u(z) = \amalg_{x\in\Crit_kf}
        n(z,x)\cdot\overline{W}{}^u(x)\,,
$$
such that by setting
$$
        U=\amalg_{\mu(x)=k} a_x\cdot W^u(x)
$$
we obtain the smooth $(k+1)$-dimensional manifold $W$ with boundary
$\partial W=U$. Note that $U\subset Z(a)$ is a $k$-dimensional open
submanifold with
$$
        E(Z(a)^\infty) \subset E(U^\infty)
$$
and, analogously to (\ref{eq cover}), $E(U^\infty)$ and $E(W^\infty)$
are covered by the at most $(k-1)$-dimensional submanifolds
$$
        E(W^\infty) \cup E(U^\infty) \subset
        \bigcup_{\mu(y)\leqslant k-1} \im E_y\,.
$$
Altogether, we have $E(W^\infty)\cap E(U)=\emptyset$ and
$H_k(E(W^\infty);\bbz)=0$ so that Theorem \ref{prop cycle1} (b) is applicable.
\end{proof}

In order to obtain a homomorphism $\Phi\colon
H^{\text{Morse}}_\ast(M)\to H^{\text{sing}}_\ast(M)$ we have to show
that the linear maps $\Phi_{f,g}$ are compatible with the canonical
isomorphisms from Theorem \ref{thm continuation}
$$
\Phi_{10}\colon H_\ast(f_0,g_0)\stackrel{\cong}{\longrightarrow}
H_\ast(f_1,g_1)
$$
from Theorem \ref{thm continuation}. For this purpose we first present 
an alternative construction of a pseudo-cycle associated to a
$(f,g)$-Morse cycle representing the same singular homology class.

Let us consider a smooth $1$-parameter family $(f_s,g_s)$ of functions
and Riemannian metrics with $-\infty<s\leqslant 0$ such that for some
$R>0$ 
$$
(f_s,g_s)=(f,g),\quad \text{for all}\; s<-R\,.
$$ 
The pair $(f,g)$ is Morse-Smale
as above. We now redefine for $x\in\Crit f$
$$
        W^u(x) = \{\,\gamma\colon (-\infty,0]\to M\,|\,
        \dot\gamma(s)+\nabla_{g_s}f_s(\gamma(s))=0,\,
        \gamma(-\infty)=x\,\}\,.
$$
All statements about the prior unstable manifolds remain valid for this
non-auto\-no\-mous flow. Observe that we have the following weak
compactness result: Given any sequence $(\gamma_n)\subset W^u(x)$ which
contains no convergent subsequence, there exists a subsequence $(n_k)$
and a reparametrization sequence $\tau_k\to-\infty$ such that
$$
        \gamma_{n_k}\stackrel{C^\infty_{\text{loc}}}{\longrightarrow}
        \gamma\in W^u(y)\quad\text{and}\quad
        \tau_k\ast\gamma_{n_k}\stackrel
        {C^\infty_{\text{loc}}}{\longrightarrow} u\in M_{x',y'}(f,g)\,.
$$
for some $x,x',y',y\in\Crit f$ with $\mu(x)\geqslant
\mu(x')>\mu(y')\geqslant\mu(y)$.  Again, in general, $\gamma_{n_k}$ can
converge weakly towards a multiply broken trajectory. If
$\mu(y)=\mu(x)-1$ we have $x'=x$ and $y'=y$.

Carrying out the same constructions as before, now based on the deformed
unstable manifolds associated to $(f_s,g_s)$, we obtain pseudo-cycles
$$
        \tilde E\colon \tilde Z(a) \to M
$$
associated to Morse-cycles $a\in Z_k(f,g)$ thus leading to a
homomorphism
\begin{equation}\label{eq defhomom}
        \widetilde\Phi_{f,g}\colon H_\ast(C_\ast(f,g),\partial) \to
        H_\ast(M;\bbz)\,.
\end{equation}
\begin{lemma}\label{prop homotopy1}
        The homomorphisms $\Phi_{f,g}$ and $\widetilde\Phi_{f,g}$ are
        identical. 
\end{lemma}
\begin{proof}
We have to show that the pseudo-cycles $E\colon Z(a)\to M$ and $\tilde
E\colon \tilde Z(a)\to M$ can be related by a suitable pseudo-cycle cobordism such
that Theorem \ref{prop cycle1} (b) applies.

Since asymptotically constant families $(f_s,g_s)_{s\in(-\infty,0]}$
as above form a convex set we
can consider the continuous path
$$
        (f_s,g_s)_\lambda = (1-\lambda)(f,g) + \lambda
        (f_s,g_s),\quad \lambda\in I=[0,1]\,.
$$
Given $x\in\Crit f$, the space
$$
        W^u_I(x) = \{\,(\lambda,\gamma)\,|\,
        \lambda\in I,\,\gamma\in W^u_\lambda(x)\,\}
$$
is a smooth manifold of dimension $\mu(x)+1$, where $W^u_\lambda(x)$ is
the unstable manifold of $x$ associated to the pair
$(f_s,g_s)_\lambda$. Its boundary is the disjoint union of $W^u_1(x)$
and $W^u_0(x)^\ast$, i.e. the latter with reversed
orientation. Analyzing the non-compactness of $W^u_I(x)$, we consider a
sequence $(\lambda_n,\gamma_n)$ which contains no convergent
subsequence. There exists a subsequence $(n_k)$ such that
$\lambda_{n_k}\to\lambda$ and either
$$
        \gamma_{n_k}\rightharpoonup (u,\gamma)\in \widehat
        M_{x,y}(f,g)\times W^u_\lambda(x)\,,
$$
for $\mu(y)=\mu(x)-1$, or $\gamma_{n_k}$ converges in
$C^\infty_{\text{loc}}$ towards a $\gamma\in W^u_\lambda(z)$ with $\mu(z)\leqslant\mu(x)-2$.

Moreover, we can prove a $\lambda$-parameterized version of the gluing
result in Lemma \ref{prop gluing1} yielding a gluing map
$$
        \#^V\colon \widehat M_{x,y}\times V\times [\rho_V,\infty)
        \to W^u_I(x)
$$
for every relatively compact, open subset $V\subset W^u_I(y)$ and
$\mu(y)=\mu(x)-1$. This allows us to build $\overline{W}_I^u(x)$ as in
(\ref{eq boundary}) and to construct a smooth manifold
$Z_I(a)$ together with a smooth map
$$
        E\colon Z_I(a)\to M,\quad
        \partial Z_I(a)=\tilde Z(a) - Z(a),
$$
which extends the given maps $\tilde E$ and $E$ on the boundary. Since
$$
        E(Z_I(a)^\infty) \subset \bigcup_{\mu(z)\leqslant\mu(x)-2}
        \im E_z
$$
for $E_z\colon W^u_I(z)\to M$ with $\dim W^u_I(z)\leqslant \mu(x)-1$, we
meet the conditions of Theorem \ref{prop cycle1} (b).
\end{proof}
In view of (\ref{eq Morse homology}), we now show that the
pseudo-cycle homomorphisms $\Phi_i=\Phi_{(f^i,g^i)}$ are compatible
with the canonical isomorphisms $\Phi_{ij}$,
\begin{lemma}\label{prop natural}
  The homomorphisms $\Phi_i\colon H_\ast(f^i,g^i) \to H_\ast(M;\bbz)$
  are compatible with $(\Phi_{ij})$, that is,
$$
        \Phi_1 \circ \Phi_{10} = \Phi_0
$$
  for all Morse-Smale pairs $(f^0,g^0)$ and $(f^1,g^1)$.
\end{lemma}
\begin{proof}
We have to compare the pseudo-cycles
$$
        E^0\colon Z^0(a) \to M \quad\text{and}\quad
        E^1\colon Z^1(\Phi_{10}(a)) \to M
$$
for any Morse cycle $a\in Z_k(f^0,g^0)$. That is, we have to show that the $E^i$
can be extended to a suitable cobordism $E\colon W\to M$ such that,
again, Theorem \ref{prop cycle1} (b) applies.

Let us consider the space similar to $W^u_I(x)$ in the proof of
Lemma \ref{prop homotopy1},
$$
        W_{\bbr_+}(x)=\{\,(\lambda,\gamma)\,|\,
        \lambda\in[0,\infty),\,\gamma\in
        W^u(x;f_{s+\lambda},g_{s+\lambda})\,\}
$$
for $x\in \Crit f^0$. (If $[0,\infty)$ is replaced by a compact interval
we are in the situation of Lemma \ref{prop homotopy1}.)
Now we have to deal with additional non-compactness for
$\lambda_n\to\infty$. Let $(\lambda_n,\gamma_n)\subset W_{\bbr_+}(x)$ be
such that $\lambda_n\to\infty$.  Then, there exists a subsequence
$(n_k)$ such that
$$
        \gamma_{n_k}\stackrel{C^\infty_{\text{loc}}}
        {\longrightarrow} \gamma\in W^u(x';f^1,g^1)
$$
for some $x'\in \Crit f^1$. Necessarily, $\mu(x')\leqslant
\mu(x)$. If both critical points $x$ and $x'$ have equal Morse index
then, up to choosing a subsequence, 
$$
        (-\lambda_{n_k})\ast\gamma_{n_k} \stackrel
        {C^\infty_{\text{loc}}}{\longrightarrow}
        u\in M_{x,x'}(f_s,g_s)\,.
$$
In that case we denote this weak convergence again
by
$$
        (\lambda_{n_k},\gamma_{n_k})\rightharpoonup
        (u,\gamma)\,.
$$
For the converse, we have a gluing theorem analogous to Lemma
\ref{prop gluing1}: 

Let $\mu(x)=\mu(x')$. Given $V\subset W^u(x')$, an open and relatively
compact subset, there exists a $\lambda_V>0$ and a smooth map
$$
        \#^V\colon M_{x,x'}(f_s,g_s)\times V\times
        [\lambda_V,\infty) \to W_{\bbr_+}(x)\,,
$$
such that the corresponding properties (a)--(d) as in Lemma
\ref{prop gluing1} hold true. Extending the construction from the proof
of Lemma \ref{prop homotopy1} based on the $\lambda$-parametrized
gluing, we now glue in boundary manifolds to $W_{\bbr_+}$ such that
\begin{align*}
        \overline{W}_{\bbr_+}(x)= W_{\bbr_+}(x) \,&\cup
        \bigcup_{\mu(y)=\mu(x)-1} \big(M_{x,y}(f^0,g^0)\times
        W_{\bbr_+}(y)\big)\\
        & \cup \bigcup_{\mu(x')=\mu(x)} \big(M_{x,x'}(f_s,g_s)\times
        W^u(x';f^1,g^1)\big)\,.
\end{align*}
Note that it is not necessary to glue in the codimension-$2$ manifolds
$M_{x,y}(f^0,g^0)\times M_{y,y'}(f_s,g_s)\times W^u(y';f^1,g^1)$ for
$\mu(y')=\mu(x)-1$. Building the quotient manifold
$$
        \bar Z(a) = \amalg_{x\in\Crit_kf^0}
        a_x\cdot \overline{W}_{\bbr_+}(x) \big/ \sim
$$
analogously as above, we obtain a smooth $(k+1)$-dimensional manifold with
boundary
$$
        \partial \bar Z(a) = Z^o(\Phi_{10}(a)) - Z(a)
$$
where $Z^o(\Phi_{10}(a))$ is the open subset
$$
        \bigcup_{x\in\Crit f^0} a_x\cdot
        M_{x,x'}\times W^u(x') \subset Z(\Phi_{10}(a))
$$
with complementary strata of codimension at least $1$. The evaluation
maps $E(\gamma)=\gamma(0)$ extend from the boundary manifolds to $\bar
Z(a)$ and it is straightforward to verify that the conditions for
Theorem  \ref{prop cycle1} (b) are satisfied.
\end{proof}
Summing up, we obtain the well-defined homomorphism
\begin{equation}\label{eq morphism}
  \Phi\colon H^{\text{Morse}}_\ast(M;\bbz)\to
  H^{\text{sing}}_\ast(M;\bbz)\,.
\end{equation}

\subsubsection{Remarks on compatibility with other equivalences for
  Morse homology}\label{ssc other}
In view of the axiomatic approach to Morse homology adopted in
\cite{Sch-Morse}, it is straightforward, based on Lemma \ref{prop
natural}, to verify that the homomorphism $\Phi$ is natural. This
means it respects functoriality with respect to closed
embeddings, w.r.t.~changes of Morse functions, and it is compatible with the
relative version of Morse homology. Thus, we can refer to the
uniqueness result from \cite{Sch-Morse}, mentioned in \ref{ssc axiom},
in order to conclude that $\Phi$ is in fact the unique, natural
isomorphism between Morse homology and singular homology.

Let us also remark that in case of a self-indexing Morse function $f$, i.e.
$$
\mu(x)=f(x), \;\forall x\in\Crit f,
$$
we have the obvious isomorphism
\begin{equation}\label{eq self-index}
\Gamma_k\colon C_k(f)\stackrel{\cong}{\longrightarrow}
H^{\text{sing}}_k(M^k,M^{k-1};\bbz),
\end{equation}
for $M^a=\{\,p\in M\,|\,f(p)\leqslant a\,\}$, $a\in\bbr$. A classical
proof for the equivalence of Morse homology and singular
homology\footnote{used in \cite{Mil-hcobordism}} is to show 
\begin{equation}\label{eq compat1}
\Gamma_{k-1}\circ \partial(f,g)=\partial_\ast\circ\Gamma_k
\end{equation}
for the boundary operator in the long exact sequence associated to the
decomposition $(M^k,M^{k-1})_{k=0,\ldots,n}$,
\begin{equation}\label{eq compat2}
 H_k(M^k,M^{k-1})\stackrel{\partial_\ast}{\longrightarrow}H_{k-1}(M^{k-1},M^{k-2})\,.
\end{equation}
Obviously, we have for $j\colon H_k(M^k)\to H_k(M^k,M^{k-1})$,
\begin{equation}\label{eq compat3}
  j\circ\Phi_{fg}(a)=\Gamma(a),\;\forall a\in Z_k(f,g)
\end{equation}
so that the induced isomorphism $\Gamma\colon
H_\ast(f,g)\stackrel{\cong}{\longrightarrow}H^{\text{sing}}_\ast(M;\bbz)$
and $\Phi_{f,g}$ are identical.

\subsection{The Inverse Homomorphism}
Although it is already clear that the homomorphism 
$$
\Phi\colon H^{\text{Morse}}_\ast(M;\bbz)\to H^{\text{sing}}_\ast(M;\bbz)
$$
is a natural isomorphism, let us nevertheless construct its inverse $\Psi=\Phi^{-1}$
{\em explicitly} along the same lines as used for $\Phi$. 

The main idea is to define an intersection number for pseudo-cycles and 
stable manifolds $W^s(x)$ for a generic Morse-Smale pair $(f,g)$. 
Recall from \cite{Sch-Morse} the construction of a Banach manifold $\Gcal$ of smooth Riemannian
metrics, $L^2$-dense in  the space of all smooth Riemannian
metrics. We consider a metric $g$ to be generic with respect to a
certain property, if we can find a residual set $\Rcal\subseteq \Gcal$ 
of metrics with that property. The first step is
\begin{theorem}\label{thm psistable}
  Let $\chi\colon V^k\to M$ be a smooth map of a $k$-dimensional
  manifold into $M$ and $f$ a Morse function on $M$ such that
  $\chi(V)\cap\Crit f=\emptyset$ if $k<n$ or $\rank D\chi(p)=n$ for
  all $p\in\chi^{-1}(\Crit f)$ if $k=n$. Then there exists a residual set
  $\Rcal\subseteq\Gcal$ such that
$$
\Mcal_{\chi;x}(f,g)=\{\,(p,\gamma)\in V\times
W^s(x)\,|\,\dot\gamma+\nabla_g f(\gamma)=0,\;
\gamma(0)=\chi(p)\,\}
$$
  is a smooth manifold of dimension
$$
\dim \Mcal_{\chi;x}(f,g)=k-\mu(x)
$$
  for all $x\in\Crit f$ and $g\in\Rcal$. In particular, it is empty if
  $k<\mu(x)$.
\end{theorem}
Also, as will be clear from the proof, if $V^k$, $M$ and $W^s(x)$ are
oriented, the intersection manifold $\Mcal_{\chi;x}(f,g)$ inherits a
well-defined orientation which is a number $\pm1\in\bbz_2$ if $k=\mu(x)$.
\begin{proof}
The main ingredient of this transversality theorem is the following
\begin{lemma}\label{lm submersion}
  The universal stable manifold
$$
W^s_{\text{univ}}\{\,(\gamma,g)\in
C^\infty([0,\infty),M)\times\Gcal\,|
\dot\gamma+\nabla_g f(\gamma)=0,\,\gamma(+\infty)=x\,\}
$$
  for $x\in\Crit f$, $f$ a Morse function, admits a submersion
$$
E\colon W^s_{\text{univ}}(x)\to M,\quad E(\gamma,g)=\gamma(0),
$$
away from the critical point $\gamma\equiv x$. It is also a
submersion everywhere if $\mu(x)=0$.
\end{lemma}
\begin{proof}
Let us recall some analytic constructions from \cite{Sch-Morse}. The
space
$$
\Hcal^{1,2}_x=H^{1,2}_x([0,\infty),M)
$$
is the $H^{1,2}$-Sobolev completion of the space of smooth curves
$\gamma\colon[0,\infty)\to M$ with sufficiently fast convergence toward $x\in M$
as $t\to\infty$. It is in fact a Hilbert manifold. The tangent space to the Banach manifold of smooth
Riemannian metrics on $M$ is
$$
T_g\Gcal=\{\,h\in C^\infty_\epsilon(\End(TM))\,|\,h\,\text{symmetric
  w.r.t. }g_o\,\}
$$
for some fixed Riemannian metric $g_o$. The function space
$C^\infty_\epsilon$ is an $L^2$-dense subspace of $C^\infty$ with a
Banach space norm. Let us now consider the smooth map
\begin{gather*}
  F\colon \Hcal^{1,2}_x\times\Gcal\to L^2(\Hcal^{1,2}_x{}^\ast TM),\\
  F(\gamma,g)=\dot\gamma+(\nabla_g f)\circ\gamma,
\end{gather*}
where the right hand side space is a Banach space bundle over the
manifold $\Hcal^{1,2}_x$ with fiber $L^2(\gamma^\ast TM)$ of
$L^2$-vector fields along the curve $\gamma$. Choosing a Riemannian
connection $\nabla$ on $TM$ we obtain the linearization of $F$ as
\begin{align*}
  DF(\gamma,g)(\xi,h)&= DF_1(\gamma,g)(\xi)+ DF_2(\gamma,g)(h),\\
  DF_1(\gamma,g)(\xi)&=\nabla_t\xi+(\nabla_\xi\nabla_g
  f)\circ\gamma,\\
  DF_2(\gamma,g)(h)&=h(\gamma)\cdot\nabla_g f(\gamma)\,.
\end{align*}
Observe that $h(\gamma)$ is an endomorphism of the pull-back bundle
$\gamma^\ast TM$. Hence, any variation of $h(\gamma(t))$ as  a
function of time $t$ can be achieved through a variation of $h$ over
$M$ if $\gamma$ is injective. Altogether we obtain the tangent space
of the universal stable manifold as
\begin{equation}\label{eq tangent}
  T_{(\gamma,g)}W^s_{\text{univ}}(x)=\{\,(\xi,h)\,|\,DF(\gamma,g)(\xi,h)=0\,\},
\end{equation}
because $0$ is a regular value for $F$, as it will be clear below.
Given $\gamma(0)=p\in M$ and $(\gamma,g)\in W^s_{\text{univ}}(x)$ we
have to show that for each $v\in T_p M$ there exist $(\xi,h)\in
T_{(\gamma,g)}W^s_{\text{univ}}(x)$ such that $\xi(0)=v$, if either
\begin{enumerate}
\item[(a)] $\gamma(0)\not\in\Crit f$, i.e.~$\gamma(0)\not=x$, or
\item[(b)] $\gamma(0)=x$ and $\mu(x)=0$.
\end{enumerate}
In the latter case (b) we have $\gamma\equiv\const=x$ with
$$
DF_1(x,g)(\xi)=\dot\xi+\Hess f(x)\cdot\xi,\quad DF_2(x,g)(h)=0,
$$
where the Hessian at $x\in\Crit_0 f$ is positive definite. This
implies
$$
\ker DF(x,g)=T_xM\times T_g\Gcal,
$$
and hence the submersion property of $E$.

In case (a) let us simplify the operator $DF(\gamma,g)$ by using
coordinates with respect to an orthonormal parallel frame of
$\gamma^\ast TM$. We obtain the operator,
\begin{equation}\label{eq parallel}\begin{split}
  D&\colon H^{1,2}([0,\infty),\bbr^n)\times T_g\Gcal\to L^2([0,\infty),\bbr^n)\\
  D&(\xi,h)=\dot\xi+A(t)\xi+h\cdot X,
\end{split}
\end{equation}
where $A\colon [0,\infty)\to S(n,\bbr)$ is a smooth path in the space
of symmetric $n\times n$-matrices with $A(\infty)=\Hess f(x)$ and
$X\colon [0,\infty)\to\bbr^n$ with $X(t)\not=0$ for all
$t\in[0,\infty)$.
We shall now prove that for all $\eta\in L^2([0,\infty),\bbr^n)$ and
$v\in\bbr^n$ there exist $\xi\in H^{1,2}([0,\infty),\bbr^n)$ and $h\in
C^\infty_o([0,\infty),S(n))$ such that $D(\xi,h)=\eta$ and $\xi(0)=v$.
This concludes the proof of (a) in view of the fact that each such $h$
arises from an $h\in T_g\Gcal$ since $\gamma$ is injective if
$\gamma(0)\not=x$.

Suppose that there exist $\eta$ and $v$ such that
\begin{equation}\label{eq cokern}
  \langle D(\xi,h),\eta\rangle_{L^2}+\langle
  \xi(0),v\rangle_{\bbr^n}=0\quad\text{for all }\xi,h\,.
\end{equation}
This implies that $\eta\in H^{1,2}([0,\infty),\bbr^n)$ and
$\dot\eta-A^t(t)\eta=0$ and therefore $\eta\equiv 0$ if $\eta(0)=0$.
Moreover, (\ref{eq cokern}) implies that $\langle hX,\eta\rangle=0$ for
all $h$ and we have $X(t)\not=0$. If $\eta(0)\not=0$ we can find\footnote{Compare (2.38)
in the proof of Proposition 2.30 in \cite{Sch-Morse}} $h(t)$
with support in $[0,\epsilon)$ such that $\langle
hX,\eta\rangle\not=0$ contradicting (\ref{eq cokern}). Hence we obtain
$\eta\equiv 0$ and by (\ref{eq cokern}) $\langle
\xi(0),v\rangle_{\bbr^n}=0$ for all $\xi$ which implies $v=0$.
Since the cokernel of $D$ in $L^2$ is finite-dimensional it follows
that $D$ is surjective.
\end{proof}

The proof of Theorem \ref{thm psistable} now follows from the parameter
version of the Sard-Smale theorem. There exists a residual subset
$\Rcal\subseteq\Gcal$ such that for each $g\in\Rcal$ the map 
$$
 (\chi,E)\colon V^k\times W^s_g(x)\to M\times M
$$
intersects the diagonal $\triangle=\{\,(p,p)\,|\,p\in M\}$
transversely. For such a generic $g$, 
$$
\Mcal_{\chi;x}(f,g)=(\chi,E)^{-1}(\triangle)
$$
is a smooth manifold of dimension $k+(n-\mu)-n$.

Note that for $k=\mu(x)$ intersections $(p,\gamma)\in \Mcal_{\chi;x}(f,g)$ for a
regular $g$ can only occur if $\rank D\chi(p)=k$. Therefore, it is
obvious how the solution space $\Mcal_{\chi;x}$ inherits its
orientation from an orientation of $V^k$, $M$ and $W^s(x)$.
\end{proof}

The main consequence of the intersection theorem \ref{thm psistable} is the compactness result
\begin{corollary}\label{cor intersect}
  For each $k$-dimensional pseudo-cycle $\chi\colon V^k\to M$
  with $\overline{\chi(V)}\cap\Crit f=\emptyset$ if $k<n$, and $\rank
  D\chi(p)=n$ for all $p\in\chi^{-1}(\Crit f)$ and $\chi(V^\infty)\cap\Crit 
  f=\emptyset$ if $k=n$, there is a
  residual set of metrics $\Rcal$ such that the intersection set
  $\Mcal_{\chi;x}(f,g)$ is finite for all 
  $x\in\Crit f$ with $\mu(x)=k$ and $g\in\Rcal$.
\end{corollary}
\begin{proof}
Consider a sequence
$(p_n,\gamma_n)\subseteq\Mcal_{\chi;x}(f,g)$. After choosing a
suitable subsequence we have
\begin{equation}\label{eq intersect1}
  \chi(p_n)\to x_o\in \overline{\chi(V)},\quad
  \gamma_n\stackrel{C^\infty_{\text{loc}}}{\longrightarrow}\gamma_o\in 
  W^s(x'),\;\mu(x')\geqslant \mu(x),\;\gamma_o(0)=x_o\,.
\end{equation}
In the case that $x_o\in \chi(V^\infty)$, we use that $\chi(V^\infty)$ 
can be covered by a map $\tilde\chi\colon \tilde V^{k-2}\to M$, so
that $(p_o,\gamma_o)\in\Mcal_{\tilde\chi;x'}(f,g)$, $\tilde\chi(p_o)=x_o$. Since the
intersection of residual sets is residual it follows from Theorem
\ref{thm psistable} that for a generic $g$ $\Mcal_{\tilde\chi;x'}$ has 
to be empty by dimensional reasons. Thus $x_o\in \chi(V^\infty)$ can
be excluded.

Sharpening the convergence result (\ref{eq intersect1}) we can deduce
weak convergence towards a broken trajectory 
\begin{equation}\label{eq intersect2}
\gamma_n\rightharpoonup (\gamma_o,u_1,\ldots,u_r),\quad
\gamma_o\in W^s(x'),\;u_1\in\Mcal_{x',x_1},\ldots,\;
u_r\in\Mcal_{x_{r-1},x}\,.
\end{equation}
For such multiply broken trajectories we must have
$\mu(x_{i-1})>\mu(x_i)$. Hence, if $k=\mu(x)$ we cannot have
$x'\not=x$ and $\Mcal_{\chi;x}(f,g)$ must be compact and hence finite.
\end{proof}

Applying the concept of coherent orientations we can now
associate to each intersection $(p,\gamma)\in\Mcal_{\chi;x}(f,g)$ a sign
$\tau(p,\gamma)$. Given a $k$-dimensional pseudo-cycle $\chi\colon
V^k\to M$ representing a singular cycle $\alpha=\alpha_\chi\in
H^{\text{sing}}_k(M)$ with $k<n$ we can find a Morse function $f$ such 
that $\Crit f\cap\overline{\chi(V)}=\emptyset$. If $k=n$, after
possibly homotoping $\chi$ to a suitable cobordant pseudo-cycle, we can 
find a Morse function $f$ such that we have only $p\in\chi^{-1}(\Crit f)$ 
with $\rank D\chi(p)=n$. We then define in view of Theorem \ref{thm
  psistable} for a generic $g$
\begin{equation}\label{eq psi}
  \Psi(\chi)=\sum_{x\in\Crit_k f}
  \#_{\text{alg}}\Mcal_{\chi;x}(f,g)\,x\,\in C_k(f,g)\,.
\end{equation}
\begin{corollary}\label{thm well-defined}
  The chain $\Psi(\chi)\in C_k(f,g)$ is a Morse-cycle, and
  given two cobordant pseudo-cycles $\chi,\chi'$, the associated Morse-cycles are
  cohomologous,
  $\Psi(\chi)-\Psi(\chi')=\partial(f,g)b$ for some
  $b\in C_{k+1}(f,g)$.
\end{corollary}
\begin{proof}
Computing
\begin{gather*}
\partial(f,g)\sum_{x\in\Crit_k f} \#_{\text{alg}}\Mcal_{\chi;x}(f,g)\,x=
\sum_{\mu(y)=k-1} n(\chi;y)\,y,\\
n(\chi;y)=\sum_{\mu(x)=k}\#_{\text{alg}}\Mcal_{\chi;x}(f,g)\#_{\text{alg}}\Mcal_{x;y}(f,g),
\end{gather*}
we have to show that
\begin{equation}\label{eq well1}
  n(\chi;y)=0\,.
\end{equation}
This follows readily from the $1$-dimensional compactness result for
$\Mcal_{\chi;y}$ analogous to (\ref{eq intersect1}) and (\ref{eq
  intersect2}) and the corresponding gluing operation. Namely, since 
$\dim\chi-\mu(y)=1$, non-compactness of $\Mcal_{\chi;y}(f,g)$ for
generic $g$ can only occur in terms of simply broken trajectories in
the limit. But exactly as for the proof of the fundamental fact
$\partial(f,g)^2=0$, the corresponding gluing result completely
analogous to Lemma \ref{prop gluing1} shows that the oriented number
of boundary components of $\Mcal_{\chi;y}(f,g)$ equals $n(\chi;y)$ and 
has to vanish since each component of $\Mcal_{\chi;y}$ is
diffeomorphic to an interval. This proves (\ref{eq well1}).

Given a pseudo-cycle cobordism $F\colon W^{k+1}\to M$ in the sense of
Theorem \ref{prop cycle1} (b), i.e. $\partial F=\chi-\chi'$, we can define 
the $1$-dimensional manifold $\Mcal_{F;x}(f,g)$ for $x\in\Crit_k f$,
and generic $g$. The same compactness-gluing argument as before shows
\begin{equation}
  \partial(f,g)\sum_{x\in\Crit_k
    f}\#_{\text{alg}}\Mcal_{F;x}(f,g)\,x=\Psi(\chi)-\Psi(\chi')\,.
\end{equation}
\end{proof}
The thus well-defined homomorphism
$$
\Psi_{f,g}\colon H^{\text{sing}}_\ast(M;\bbz)\to H_\ast(f,g)
$$
is compatible with the canonical isomorphism
$$
\Phi_{10}\colon H_\ast(f^0,g^0)\to H_\ast(f^1,g^1)\,.
$$
We have
\begin{corollary}\label{cor isomorphism}
  Considering the isomorphism $\Phi_{10}$ for generic Morse-Smale
  pairs, it holds
$$
\Phi_{10}\circ\Psi_{f^0,g^0}=\Psi_{f^1,g^1},
$$
  and we have for the well-defined homomorphism 
$$
\Psi\colon H^{\text{sing}}_\ast(M;\bbz)\to
H^{\text{Morse}}_\ast(M;\bbz)
$$
  the identity $\Psi\circ\Phi=\id_{H^{\text{Morse}}}$.
\end{corollary}
\begin{proof}
The proof of the compatibility
$\Phi_{10}\circ\Psi_{f^0,g^0}=\Psi_{f^1,g^1}$ can be carried out
exactly analogous to that for Lemma \ref{prop natural} using the
argument from the proofs of Corollaries \ref{cor intersect} and
\ref{thm well-defined}.

Consider now the $k$-dimensional pseudo-cycle
$E\colon Z(a)\to M$ associated to a Morse-cycle $a\in Z_k(f^0,g^0)$. Let
$y\in\Crit_k f^1$. Then in view of Theorem \ref{thm psistable} for a
Morse-Smale pair $(f^1,g^1)$ with generic $g^1$ we have intersections
$(p,\gamma)\in\Mcal_{E;y}(f^1,g^1)$ only for $p\in Z(a)$ on the
$k$-dimensional strata which are exactly the unstable manifolds
$W^u(x,f^0,g^0)$ in $a=\sum_{x\in\Crit_k f^0}a_x\,x$. We therefore
have 
$$
\Mcal_{E;y}(f^1,g^1)=\{\,(\gamma^-,\gamma^+)\in W^u(x,f^0,g^0)\times
W^s(y;f^1,g^1)\,|\, \gamma^-(0)=\gamma^+(0)\,\}\,.
$$
Using a homotopy operator as before we can show easily that this 
is homologically equivalent to the definition of the operator $\Phi_{10}\colon
C_\ast(f^0,g^0)\to C_\ast(f^1,g^1)$, i.e.
\begin{equation}\label{eq iso1}
  \Psi_{f^1,g^1}\circ\Phi_{f^0,g^0}=\Phi_{10}\colon
  H_\ast(f^0,g^0)\stackrel{\cong}{\longrightarrow}H_\ast(f^1,g^1)\,.
\end{equation}
This proves $\Psi\circ\Phi=\id_{H^{\text{Morse}}}$.
\end{proof}
Using the fact that $H^{\text{Morse}}_\ast(M)\cong
H^{\text{sing}}_\ast(M)$, it follows immediately that the left-inverse 
$\Psi$ is the inverse of $\Phi$. Thus we have explicit constructions
of both isomorphisms in terms of Morse-pseudo-cycle equivalences.
\bibliographystyle{amsalpha}
\providecommand{\bysame}{\leavevmode\hbox to3em{\hrulefill}\thinspace}

\end{document}